\documentclass{elsart}
\usepackage{amssymb, amsmath}
\newcommand{\be}{\begin{equation}}
\newcommand{\ee}{\end{equation}}
\newcommand{\xx}{{\mathbf x}}
\renewcommand{\aa}{{\mathbf a}}
\newcommand{\bb}{{\mathbf b}}
\newcommand{\sumi}{{\sum_{i=1}^n}}
\begin{document}

\begin{frontmatter}

\title{Inference for Multivariate Normal Mixtures}

\author{Jiahua Chen}
\address{Department of Statistics, University of British Columbia \\ Vancouver, BC, V6T 1Z2, Canada}

\ead{jhchen@stat.ubc.ca}

\author{Xianming Tan}
\address{LMPC and School of Mathematical Sciences, Nankai University \\ Tianjin, 300071, P.R. China} %

\ead{tanxm@nankai.edu.cn}

\begin{abstract}
Multivariate normal mixtures provide a flexible model for
high-dimensional data. They are widely used in statistical genetics,
statistical finance, and other disciplines. Due to the unboundedness of
the likelihood function, classical likelihood-based methods, which
may have nice practical properties, are inconsistent. In this paper,
we recommend a penalized likelihood method for estimating the mixing
distribution. We show that the maximum penalized likelihood
estimator is strongly consistent when the number of components has a
known upper bound. We also explore a convenient EM-algorithm for
computing the maximum penalized likelihood estimator. Extensive
simulations are conducted to explore the effectiveness and the
practical limitations of both the new method and the ratified maximum
likelihood estimators. Guidelines are provided based on the simulation
results. %
\end{abstract}

\begin{keyword}
Multivariate normal mixture \sep Penalized maximum likelihood
estimator \sep Strong consistency.

\PACS 02.50.-r
\end{keyword}
\end{frontmatter}



\section{Introduction}

In the past few decades, there has been an exploding volume of
literature on mixture models [\ref{titt}, \ref{lind1}, \ref{mac}, \ref{fru}].
Various mixture distributions including normal mixtures
are used in a wide variety of situations. %
Schork et al.\ [\ref{sch}] reviewed the applications of mixture
models in human genetics and Tadesse et al.\ [\ref{tad}] used a normal mixture
model for clustering analysis. Application examples can be found in
[\ref{fra}, \ref{lin}, \ref{raf}]
and [\ref{alex}]. 

Finite mixtures of multivariate normals have also drawn substantial
attention recently. Lindsay and Basak [\ref{lind2}] devised a system of
moment equations and a fast algorithm to estimate the parameters of
multivariate normal mixture distributions under an 
equal-covariance-matrix assumption. 
However the equality assumption is crucial, and
failing this condition leads to a substantial loss in the accuracy of
the fit [\ref{mac}]. Unequal-variance normal
mixture models have an ill effect on the likelihood function
[\ref{day}]. Placing a positive lower bound on the component variances helps,
but the resulting statistical procedure can be awkward because it is
not continuous in the data. Placing a positive lower bound on the ratio
of the component variances is better. In the univariate case the
resulting constrained maximum likelihood estimator is consistent for
both constant and shrinking lower bounds [\ref{hat}, \ref{tan}]. 
Though consistency is yet to be proved, Ingrassia [\ref{ing}]
applied the constrained method to multivariate observations. Ray and
Lindsay [\ref{ray}] found that in contrast to the univariate case, the multivariate
normal mixture density can have more modes than the number of
components. Inference on multivariate normal mixture models is
hence more difficult.

In this paper, we investigate a penalized likelihood method for
estimating the mixing distribution. 
The penalized likelihood estimations form a population class
of methods, see [\ref{gre}, \ref{egg}].
When the number of components
has a known upper bound, the maximum penalized likelihood estimator
(PMLE) is found to be strongly consistent. An EM-algorithm is developed
and extensive simulations are conducted. Although after some
ratification, the usual maximum likelihood estimators and the PMLE
work similarly after the removal of degenerating local maxima in the
univariate case [\ref{chen}], the PMLE is
advantageous for multivariate normal mixture models.

The paper is organized as follows. In Section 2, the penalized
likelihood method is introduced. Two theorems on strong
consistency are presented with the proofs deferred to the Appendix.
The EM-algorithm for solving the maximization problem for the penalized
likelihood function is given. Section 3 contains the simulation results.

\section{Penalized likelihood method}

\subsection{Consistency of the PMLE}

Let $\varphi({ \xx;  \mu, \Sigma})$ be the
multivariate normal density with
$(d \times 1)$ mean vector ${ \mu}$ and
$d \times d$ covariance matrix ${\Sigma}$, i.e.,
\[
\varphi({\xx;  \mu, \Sigma}) %
=%
\{ 2 \pi |\Sigma| \}^{-d/2} \exp \{ - \frac{1}{2} (\xx - \mu)^\tau \Sigma^{-1} (\xx - \mu) \}.%
\]
A $d$-dimensional random vector $X$ has a multivariate
finite normal mixture distribution of order $p$
if its density function is given by
\be%
\label{eqn1} f(\xx; G)=
\pi_1 \varphi({\xx;  \mu_1, \Sigma_1}) +%
\pi_2 \varphi({\xx;  \mu_2, \Sigma_2}) + \cdots %
\pi_p \varphi({\xx;  \mu_p, \Sigma_p})
 \ee%
where
$G$ is the mixing distribution
assigning probability $\pi_j$ to parameter set
$(\mu_j, \Sigma_j)$ of the $j$th kernel density
$\varphi(\xx; \mu_j, \Sigma_j)$.

Let $\xx_1, \xx_2, \ldots, \xx_n$ be a random sample from
(\ref{eqn1}). Then
\begin{equation*}
l_n(G) = \sumi \log f(\xx_i, G)
\end{equation*}
is the log-likelihood function.
Even if $|\Sigma_{j}| >0$ for all $j$, $l_n(G)$ is
unbounded at $\mu_1 =
\xx_1$ when $|\Sigma_1|$ gets arbitrarily small.
The penalized
log-likelihood function is of the form
\[
pl_n(G) = l_n(G) + p_n(G)
\]
where $p_n(G)$ is the penalty depending on the mixing distribution
$G$ and the sample size $n$.
Let $\hat G_n$ be the mixing
distribution in the parameter space at which $pl_n(G)$ attains its
maximum. We call $\hat G_n$ the penalized maximum likelihood
estimator (PMLE).

We choose a penalty function such that:
\begin{enumerate}
\item[C1.]
$ p_n(G)
=
\sum_{j=1}^p \widetilde{p}_n({\mathbf \Sigma_j})$, %

\item[C2.]
At any fixed $G$
such that $|\Sigma_j| > 0$ for all
$j=1, 2, \ldots, p$, we have $p_n(G) = o(n)$,
and
 $\sup_G \max \{0, p_n(G)\} = o(n)$.

In addition,
$~p_n(G)$ is differentiable with respect to $G$ and as $n \to
\infty$, $ p'_n(G)=o(\sqrt{n})$ at any fixed $G$ such that $|\Sigma_j|
> 0$ for all
$j=1, 2, \ldots, p$. Here we treat $G$ as a vector of parameters contained
in the mixing distribution $G$.%

\item[C3.]  For
large enough $n$, %
$\widetilde{p}_n(\Sigma) \leq 4(\log n)^2 \log |\Sigma|$, %
when $|\Sigma|$ is smaller than $c n^{-2d}$
for some $c>0$.

\end{enumerate}

These conditions are quite flexible and functions satisfying these
conditions can be easily constructed. A class of such functions will
be given in the simulation section. Condition C1 simplifies the
numerical computation. Condition C2 limits the effect of the penalty.
The key condition is C3: it counters the damaging effect of a
degenerate component covariance matrix. The order of the penalty
size is well calibrated as will be seen in the proof, yet the exact
value of the constant $4$ is not important. The penalty function can
also be viewed as a prior function via Bayesian analysis.

\begin{thm}
\label{Thm1} Assume that the true density function
\[
f(x; G_0) =  %
\sum_{j=1}^{p_0} \pi_{0j} \varphi(x; \mu_{0j}, \Sigma_{0j})%
\]
satisfies $\pi_{0j} > 0$, $|\Sigma_{0j}| > 0$, and $(\mu_{0j},
\Sigma_{0j}) \neq (\mu_{0k}, \Sigma_{0k})$ for all $j=1, 2,
\ldots, p_0$ and $j \neq k$.

Assume that the penalty function $p_n(G)$ satisfies C1-C3 and $\tilde
G_n$ is a mixing distribution of order $p_0$ satisfying
\[
pl_n(\tilde G_n) - pl_n(G_0) \geq c > -\infty,
\]
for all $n$. Then, as $n \to \infty$,
\(
\tilde G_n {\to} G_0, %
\) almost surely. %
\end{thm}
The proof is deferred to the Appendix.

Since $pl_n(\hat G_n) - pl_n(G_0) \geq 0$,
the PMLE $\hat G$ is strongly consistent.
Because $\hat G_n$
and $G_0$ have the same order, all elements
in $\hat G_n$ converge to those of $G_0$ almost surely.
Furthermore, let
\[
S_n(G)
= \sumi \frac{\partial \log f(x_i; G)}{\partial G}
\]
be the vector score function at $G$. Let
\[
S_n'(G) = \sumi \frac{\partial S_n(G)}{\partial G}
\]
be the matrix of the second derivative of the log-likelihood
function.
At $G = G_0$, the normal mixture model is regular and hence
the Fisher information
\[
I_n(G_0) = n I(G_0) = - E \{ S_n'(G_0)\}
= E \big [ \{ S_n(G_0) \}^\tau S_n(G_0) \big ]
\]
is positive definite. Using classical asymptotic techniques as in
[\ref{leh}], and under condition C2 such that $p'_n(G) = o_p(
n^{1/2})$, we have
\[
\hat G_n - G_0
=
\{ S_n'(G_0)\}^{-1} S_n(G_0) + o_p(n^{-1/2}).
\]
Therefore, $\hat G_n$ is an asymptotically normal and efficient
estimator.

\begin{thm}
\label{Thm2}
Under the same conditions as in Theorem \ref{Thm1},
as $n \to \infty$,
\[
\sqrt{n} \{ \hat G_n - G_0\}
\to
N(0, I(G_0))
\]
in distribution. %
\end{thm}

The proof is straightforward and omitted. In practice, we may
know only an upper bound for $p_0$ rather than its exact value.
The following theorem deals with this situation.

\begin{thm}
Assume the same conditions as in Theorem \ref{Thm1}, except
that the order of the finite normal mixture model $p_0$ is
known only to be smaller than or equal to $p$.
Let  $\tilde G_n$ be a mixing distribution of order $p$ satisfying
\[
pl_n(\tilde G_n) - pl_n(G_0) \geq c > -\infty
\]
for all $n$. Then, as $n \to \infty$, \(G_n \stackrel{w}{\to} G_0 \)
almost surely. %
\end{thm}
The proof is deferred to the Appendix.

\subsection{The EM-algorithm}
We recommend the EM-algorithm due to its simplicity in coding, and its
guaranteed convergence to some local maximum under very general
conditions [\ref{wu}, \ref{rich}, \ref{gre}].
In our simulations, we use a
number of initial values to reduce the risk of poor local maxima. We
also recommend some convenient and effective penalty functions for
the EM-algorithm.

Let $z_{ij}$ be the membership
indicator variable such that it equals 1 when
$\xx_i$ is from the $j$th component of the normal
mixture model, and equals 0 otherwise.
The complete observation
log-likelihood under a normal mixture model is
then given by
\[
l_c(G)
=
\sum_{i=1}^n \sum_{k=1}^p z_{ik}
\left \{\log \pi_k -
\frac{1}{2}\log |\Sigma_k| - \frac{1}{2} (\xx_i - \mu_k)^{\tau}
\Sigma_k^{-1}(\xx_i - \mu_k) \right \}.
\]
Given the current mixing distribution
$$
G^{(m)}
=
(\pi_1^{(m)}, \ldots, \pi_p^{(m)}, \mu_1^{(m)}, \ldots,
\mu_p^{(m)},\Sigma_1^{(m)}, \ldots, \Sigma_p^{(m)}),
$$
the EM-algorithm iterates as follows:

In the E-Step, we compute
$$
\pi_{ij}^{(m+1)}
=
E\{z_{ij}| {\xx_1, \ldots, \xx_n}, G^{(m)} \}
=
\frac{\pi_j^{(m)} \phi(\xx_i;
\mu_j^{(m)},\Sigma_j^{(m)})}
{\sum_{j=1}^p \pi_j^{(m)}
\phi(\xx_i;\mu_j^{(m)},\Sigma_j^{(m)}) }.
$$
Replacing $z_{ij}$ by $\pi_{ij}^{(m+1)}$ in $l_c(G)$, we get
\begin{eqnarray*}
Q(G; G^{(m)})
&=&
E \{ l_c(G)+ p_n(G)| \xx_1, \ldots, \xx_n, G^{(m)} \} \\
&=&
\sum_{j=1}^p  (\log \pi_j) \sum_{i=1}^n \pi_{ij}^{(m+1)}
     - \frac{1}{2} \sum_{j=1}^p \left( \log |\Sigma_j| \right) \sum_{i=1}^n \pi_{ij}^{(m+1)} \\
& &
  - \frac{1}{2} \sum_{j=1}^p
  \sum_{i=1}^n \pi_{ij}^{(m+1)} (\xx_i -\mu_j)^{\tau} \Sigma_j^{-1} (\xx_i -\mu_j)%
     +p_n(G).
\end{eqnarray*}
This completes the E-step.

In the M-step, we maximize $Q(G; G^{(m)})$ with respect to $G$ to
obtain $G^{(m+1)}$. We suggest the following penalty
functions in practice:
\begin{equation}
\label{eqn41}
 p_n(G)
 =
 -a_n  \sum_{j=1}^p \left\{ \mbox{tr}(S_x \Sigma_j^{-1})
 +
 \log |\Sigma_j| \right\}   %
\end{equation}
with $S_x$ being the sample covariance matrix, and 
$\mbox{tr$(\cdot)$ }$ being the trace function. Using this penalty function,
$Q(G; G^{(m)})$ is maximized at $G= G^{(m+1)}$ with
\begin{equation*}
\left\{ \begin{aligned}%
\pi_j^{(m+1)} %
&= %
\frac{1}{n} \sum_{i=1}^n \pi_{ij}^{(m+1)}, \\
\mu_j^{(m+1)} %
&=%
\frac{\sum_{i=1}^n \pi_{ij}^{(m+1)} \xx_i}{n \pi_j^{(m+1)}},\\ %
\Sigma_j^{(m+1)}
&=%
\frac{2 a_n S_x + S_j^{(m+1)}}{2 a_n + n \pi_j^{(m+1)}}
\end{aligned} \right.
\end{equation*}
where %
$$
S_j^{(m+1)}%
=%
\sum_{i=1}^n \pi_{ij}^{(m+1)}%
(\xx_i - \mu_j^{(m+1)})(\xx_i - \mu_j^{(m+1)})^{\tau}.
$$
From a Bayesian point of view, the penalty function
(\ref{eqn41}) puts a Wishart distribution prior on $\Sigma_j$,
and $S_x$ is the mode of the prior distribution. Increasing the
value of $a_n$ implies a stronger conviction on $S_x$ as the
possible value of $\Sigma_j$.

The EM-algorithm iterates between the E-step and the M-step. The
penalized likelihood increases after each iteration. At the same
time, the {\it penalized} likelihood is bounded over the parameter
space. Hence, the EM-algorithm converges to a non-degenerate local
maximum. This is the dividing line between the penalized
likelihood and the ordinary likelihood. In both cases,
the EM-algorithm may converge to an undesired local maxima starting
from a poor initial value. In the simulations, we use ten initial
values including the true value for each data set to 
control this potential problem.

\section{Simulation study.}
When computing the MLE the local maxima located by the EM-algorithm
with degenerate covariance matrices are first removed.  The one
that attains the largest likelihood value among those remaining
is then identified as the MLE or the ratified MLE of the mixing
distribution. Although this approach lacks solid theoretical support,
it works well for univariate normal mixture models [\ref{chen}]. 
The consistency result for the PMLE for multivariate normal
mixture models does not guarantee its superiority in practice.
Thus, we feel obliged to compare the performance of the
PMLE with that of the ratified MLE. In addition, there is a general shortage
of thorough simulation studies in the context of multivariate normal
mixture models. This paper partially fills that knowledge gap.

We use bias and standard deviation to measure the accuracy of
the ratified MLE
and the PMLE. We also record the number of times that the
EM-algorithm degenerates when the ratified MLE is attempted. For
clarity, the simulation results are organized into two subsections.

\subsection{Simulation models and settings}
The size of the parameter space for the finite multivariate normal
mixture model explodes with the dimension. It is difficult to use
a few typical specific distributions to cover all aspects of this
model. We struggled to come up with a few particularly important cases.
We considered four
categories of mixture models: two-component bivariate normal
mixture models ($p=2, \,\, d=2$); three-component bivariate normal
mixture models ($p=3, \,\, d=2$); two-component trivariate normal
mixture models ($p=2, \,\, d=3$); and three-component trivariate
normal mixture models ($p=3, \,\, d=3$).

In each category, we chose $3 \times 6 $ models
formed by component mean vector and
covariance matrix configurations.
These combinations mimic
practical situations and make the comparison of the
performance of the ratified MLE and the PMLE meaningful.

The covariance matrices in the simulation models are designed
to have the following general form when $d=2$:
\[
\Sigma
=
\left [%
\begin{array}{cc}
\cos \theta & -\sin \theta\\
\sin \theta &  \cos \theta%
\end{array}
\right ]
\left [%
\begin{array}{cc}
\lambda_1 & 0\\
0 &  \lambda_2 %
\end{array}
\right ]
\left[%
\begin{array}{cc}
\cos \theta & \sin \theta\\
- \sin \theta &  \cos \theta%
\end{array}
\right ].
\]
By the choices of the eigenvalues $\lambda_1, \lambda_2$,
and the orientation angle $\theta$, we obtain various
configurations of  bivariate normal mixture models.

The covariance matrices in the simulation models are designed
to have the following general form when $d=3$:
\[
\Sigma
= P(\alpha,\beta,\gamma)
\mbox{diag} [\lambda_1, \lambda_2, \lambda_3]
P^T (\alpha,\beta,\gamma)
\]
with\vspace*{-4ex}
\begin{eqnarray*}
&&P(\alpha,\beta,\gamma) = \\
&&
\begin{bmatrix}
 \cos\alpha \cos\gamma - \cos\beta \sin\alpha \sin\gamma &
-\cos\beta \cos\gamma \sin\alpha - \cos\alpha \sin\gamma &
 \sin\alpha \sin\beta
\\
 \cos\gamma \sin\alpha + \cos\alpha \cos\beta \sin\gamma &
 \cos\alpha \cos\beta \cos\gamma - \sin\alpha \sin\gamma &
-\cos\alpha \sin\beta
\\
 \sin\beta \sin\gamma &
 \cos\gamma \sin\beta &
 \cos\beta
\end{bmatrix},
\end{eqnarray*}
{that is, a $3 \times 3$ rotation matrix.}
For each
multivariate normal mixture model, we specify the mixing
proportion, covariance matrix, and mean vector for each component.

\vspace{2ex}
\noindent{\bf Two-component bivariate normal mixture models.} %
We set the component proportions
$(\pi_1, \pi_2) = (0.3, 0.7)$. No other cases are considered.

Due to the invariance property of the multivariate normal
distribution, the distance between the
two mean vectors is the only configuration that
can make a difference. Thus, we simulated only
three pairs of mean vectors representing the situation
where two component mean vectors are in near, moderate, and
distant locations as in the following table:
\[
\begin{tabular}{|c|c|c|c|}
\hline
            & near & moderate & distant  \\
\hline
Component 1 & (0, -1)  & (0, -3)  & (0, -5) \\
Component 2 & (0, ~1) & (0, ~3) & (0, ~5) \\
\hline
\end{tabular}
\]

There are many features in the pair of
covariance matrices that may have an effect
on the performance of the ratified MLE or PMLE.
The sizes of the eigenvalues are most
important in their ratio $\lambda_2/\lambda_1$.
The angle $\theta$
determines the relative orientation between two
component densities. Our choices based on these
considerations are given in the following table:
\[
\begin{tabular}{|c|ccc|ccc|}
\hline
 & \multicolumn{3}{|c|}{Component 1} &  \multicolumn{3}{c|}{Component 2}\\
 & $\lambda_1$ & $\lambda_2$ & $\theta$
          & $\lambda_1$ & $\lambda_2$ & $\theta $\\
          \hline
    1 & 1 & 1 & $0$      & 1 & 1 & 0\\
    2 & 1 & 5 & $0$      & 1 & 1 & 0\\
    3 & 1 & 5 & $\pi/4$ & 1 & 1 & 0\\
    4 & 1 & 5 & $\pi/2$ & 1 & 1 & 0\\
    5 & 1 & 5 & $\pi/4$ & 1 & 5 & 0\\
    6 & 1 & 5 & $\pi/2$ & 1 & 5 & 0\\
\hline
\end{tabular}
\]

\vspace{2ex}
\noindent{\bf Three-component bivariate normal mixture models.} %
We set the component proportions
$(\pi_1, \pi_2, \pi_3) = (.15, .35, .50)$.
The three mean vectors
may form a straight line, an acute triangle, or an obtuse
triangle. We select three
representative ones as follows:
\[
\begin{tabular}{|c|c|c|c|}
\hline
            & straight & acute & obtuse  \\
\hline
Component 1 & (0, -2)  & (0, -2)  & (0, -2) \\
Component 2 & (0, ~0) & (3, ~0) & (1, ~0) \\
Component 3 & (0, ~2) & (0, ~2) & (0, ~2) \\
\hline
\end{tabular}
\]
We select six triplets of covariance matrices as follows:
\[
\begin{tabular}{|l|ccc|ccc|ccc|}
\hline
 & \multicolumn{3}{|c|}{Component 1} &  \multicolumn{3}{c|}{Component 2}
  &\multicolumn{3}{c|}{Component 3}\\
& $\lambda_1$ & $\lambda_2$ & $\theta$
          & $\lambda_1$ & $\lambda_2$ & $\theta$
          & $\lambda_1$ & $\lambda_2$ & $\theta $\\
    1 & 1& 1& 0 & 1 & 1 & $0$      & 1 & 1 & 0\\
    2 & 1& 1& 0 & 1 & 1 & $0$      & 1 & 5 & 0\\
    3 & 1& 1& 0 & 1 & 5 & 0& 1 & 5 & $\pi/4$\\
    4 & 1& 1& 0 & 1 & 5 & 0& 1 & 5 & $\pi/2$\\
    5 & 1& 5& 0 & 1 & 5 & $\pi/4$& 1 & 5 & $-\pi/4$\\
    6 & 1& 5& 0 & 1 & 5 & $\pi/4$& 1 & 5 & $-\pi/2$\\
  \hline
\end{tabular}
\]

\vspace{2ex}

\noindent{\bf Two-component trivariate normal mixture models.} %
We again let $(\pi_1, \pi_2) = (0.3, 0.7)$. At the same time, only
the distance between the two mean vectors matters. The two
mean vectors are chosen to be:
\[
\begin{tabular}{|c|c|c|c|}
\hline
            & near & moderate & distant  \\
\hline
Component 1 & (0, 0, -1)  & (0, 0, -3)  & (0, 0, -5) \\
Component 2 & (0, 0, ~1) & (0, 0, ~3) & (0, 0, ~5) \\
\hline
\end{tabular}
\]
The covariance matrix pairs are chosen as follows:
\[
\begin{tabular}{|l|cc|cc|}
\hline
 & \multicolumn{2}{|c|}{Component 1} &  \multicolumn{2}{c|}{Component 2}\\
& ($\lambda_1, \lambda_2, \lambda_3$) & ($\alpha, \beta, \gamma$)
           & ($\lambda_1, \lambda_2, \lambda_3$) & ($\alpha, \beta, \gamma$) \\
          \hline
    1 & (1, 1, 1)& (0, 0, 0) & (1, 1, 1) & (0, 0, 0) \\
    2 & (1, 1, 1)& (0, 0, 0) & (1, 3, 10) & (0, 0, 0) \\
    3 & (1, 3, 10)& (0, 0, 0) & (1, 3, 10) & (0, 0, 0) \\
    4 & (1, 3, 10)& (0, 0, 0) & (1, 3, 10) & ($-\pi, \pi, \pi$)/3 \\
    5 & (1, 3, 10)& (0, 0, 0) & (1, 3, 10) & ($\pi, -\pi, \pi$)/3 \\
    6 & (1, 3, 10)& (0, 0, 0) & (1, 3, 10) & ($\pi, \pi, -\pi$)/3 \\
\hline
\end{tabular}
\]

\vspace{2ex}
\noindent{\bf Three-component trivariate normal mixture models.} %
We let the component proportions $(\pi_1, \pi_2, \pi_3)$ be
$(.15, .35, .50)$. Recall that any three points fall into one plane.
Thus, the invariance property of the normal distribution allows us
to set the first entry of the mean vector to 0:
\[
\begin{tabular}{|c|c|c|c|}
\hline
            & straight & acute & obtuse  \\
\hline
Component 1 & (0, 0, -2)  & (0, 0, -2)  & (0, 0, -2) \\
Component 2 & (0, 0, ~0) & (0, 3, ~0) & (0, 1, ~0) \\
Component 3 & (0, 0, ~2) & (0, 0, ~2) & (0, 0, ~2) \\
\hline
\end{tabular}
\]
The covariance matrix triplets are chosen as follows:
\[
\begin{tabular}{|c|cc|cc|cc|}
\hline
 & \multicolumn{2}{|c|}{Component 1} &  \multicolumn{2}{c|}{Component 2}
    &  \multicolumn{2}{c|}{Component 3}\\
& ($\lambda_1, \lambda_2, \lambda_3$) & ($\alpha, \beta, \gamma$)
           & ($\lambda_1, \lambda_2, \lambda_3$) & ($\alpha, \beta, \gamma$)
           & ($\lambda_1, \lambda_2, \lambda_3$) & ($\alpha, \beta, \gamma$) \\
          \hline
    1 & (1, 1, 1)& (0, 0, 0) & (1, 1, 1) & (0, 0, 0) & (1, 1, 1) & (0, 0, 0) \\
    2 & (1, 1, 1)& (0, 0, 0) & (1, 1, 1) & (0, 0, 0) & (1, 3, 10)& (0, 0, 0) \\
    3 & (1, 1, 1)& (0, 0, 0) & (1, 3, 10) & (0, 0, 0)
                          & (1, 3, 10) & ($- \pi, \pi, \pi$)/3 \\
    4 & (1, 1, 1)& (0, 0, 0) & (1, 3, 10) & (0, 0, 0)
                          & (1, 3, 10) & ($\pi, -\pi, \pi$)/3 \\
    5 & (1, 3, 10)&(0, 0, 0)& (1, 3, 10) & ($- \pi, \pi, \pi$)/3
                          & (1, 3, 10) & ($\pi, -\pi, \pi$)/3 \\
    6 & (1, 3, 10)& (0, 0, 0)& (1, 3, 10) & ($\pi, -\pi, \pi$)/3
                          & (1, 3, 10) & ($\pi, \pi, -\pi$)/3 \\
\hline
\end{tabular}
\]

{We let $n=200$ for the two-component bivariate mixtures
and $n=300$ for the other mixtures} to ensure a reasonable estimation of
the mixing distribution. We generate 1000 data sets for each
model.

We have presented four categories of finite normal mixture models.
For ease of reference we use, for example, I.1.2 to
refer to the model from Category I with mean vector
configuration 1 and covariance
matrix configuration 2.
Even though there are many more mixing distribution
configurations for
which simulation
studies are needed, there is a limit to how much one paper can
achieve.
We do not consider the case
where $p$ is unknown.
All estimators in this case are expected to be poor
although the consistency result for the PMLE remains true.

{\bf Penalty term and initial values.}
We compute the ratified MLE and two penalized MLEs
corresponding to $a_n=n^{-1}$ and $a_n=n^{-1/2}$
in (\ref{eqn41}). We call these MLE, PMLE1, and PMLE2,
respectively.

The ten initial values are chosen from two groups. The first group
of initial values includes the true mixing distribution and four
others obtained by perturbing the
component mean vectors of the true mixing distribution. The second
group of initial values was data-based. We first calculate the
sample mean vector and the sample covariance matrix. Then we set the
mixing proportions all equal to $1/p$ and the component covariance
matrices all equal to the sample covariance matrix. We then
apply
similar perturbation
to the sample mean vector to obtain another five sets of
initial values.

\subsection{Simulation results}
{\bf Number of Degeneracies.}
When the EM-algorithm converges to a
mixing distribution with singular component
covariance matrices, we say that it degenerates.
The EM-algorithm for the PMLE
does not degenerate which is theoretically ensured.
Regardless of the quality of the initial value,
the corresponding
EM-algorithm always converges to some non-degenerate
local maximum. The PMLE is a good estimator
if the largest local maximum is a good estimator.

When computing the ratified MLE, the EM-algorithm sometimes converges
to a degenerate local maximum. We recorded the number of times
that the EM-algorithm degenerated while computing
the ratified MLE in our simulation. Since each data set had
ten initial values, the number of degenerate outcomes is out of
10,000 for each entry.

For two-component bivariate normal mixture models, it is
immediately clear that the number of degenerate outcomes increases
when the mean vectors are more widely separated. The covariance structure is
also important. For example, when the eigenvectors of one
covariance matrix are rotated by an angle of $\pi/2$ (variance
configurations 4 and 6), so that the two clusters of observations
become more mixed, the number of degenerate outcomes declines.
This observation is somewhat counter-intuitive but can be
explained as follows. The success of the EM-algorithm is heavily
dependent on sensible initial values. When the two mean vectors are
close and the components are well mixed,
different initial values do not matter as much. However, when
the two mean vectors are distant, the location of the initial mean
vectors is crucial. Thus the degenerate outcomes were mostly due to the
second
group of initial values. %

In the other three categories, the above phenomenon persists.
That is, the frequency of degeneracy increases when components are
 more widely separated. In addition, for these categories we
observe a higher frequency of degeneracies on average. We believe
this is because the EM-algorithm is more sensitive to the quality
of the initial values when the mixture models are more complicated.

{Degeneracy of the EM-algorithm should not be a}
serious problem for the ratified MLE, as long as the non-degenerate
outcomes of the algorithm provide good estimates. We hence proceed
to examine the bias and variance properties of the PMLE and the
largest non-degenerate local maxima regarded as the ratified MLE.

{\bf Bias and Standard Deviation.} We compute the element-wise
mean bias and standard deviation based on 1000 simulated samples
from each model. We
present only a subset of representative outcomes from each category; the
complete set is available upon request.

Two representative outcomes for models I.1.1 and I.2.4 in Category I
are given in Table \ref{tab4.1}. There is about a 10\% reduction in the
standard deviation for PMLE2 compared to the
ratified MLE or PMLE1 for the parameters in component 1 of Model I.1.1. The same
is true for Models I.1.5 and I.1.6 (not presented). The PMLE2
also has a relatively lower bias in these models. The results
for the remaining models are comparable to those for I.2.4: there is little
appreciable difference between the three estimation methods.

The biases of all three estimators for estimating $\mu_2$ are high
under I.1.1 and I.1.5 in which the two mean vectors are lined up in
the $\mu_1$ direction. Due to the orientation of the two component
covariance matrices, it is hard to tell the two mean vectors apart.
The biases and standard deviations for estimating $\sigma_{22}$
under I.1.1, I.1.2, $\ldots$, I.1.6 are also high or relatively high.

\vspace{3ex}%
\centerline{Table \ref{tab4.1} about here.}%
\vspace{2ex}%

We present outcomes for two models (II.1.1, II.2.4) in Category II in
Tables \ref{tab4.2a} and \ref{tab4.2b}. For both models, 
for the parameters in component
1, there is a 10\% to 20\% reduction in the standard
deviation for PMLE2 compared to the other two
estimators. The bias of PMLE2 is also lower. Some reductions in
components 2 and 3 are also noticed but to varying degrees. In the
other models, the performance of PMLE2 does not dominate that of the
ratified MLE or PMLE1.

Under a straight-line configuration of the
component mean vectors, the
bias for estimating $\mu_2$ is relatively high.
For a triangle configuration, the roles of $\mu_1$ and
$\mu_2$ are no longer different.
This bias problem is not estimator
dependent although PMLE2 helps slightly.

The estimation of $\sigma_{22}$ again comes with
both higher bias and higher standard deviation in general.
For this category of models, the problem spreads
into other parts of the covariance matrix.

\vspace{3ex}%
\centerline{Tables \ref{tab4.2a}, \ref{tab4.2b} about here.}%
\vspace{1ex}

We report simulation results for three models (III.1.1, III.2.4,
III.3.6) in Category III in Tables \ref{tab4.3a}, \ref{tab4.3b},
and \ref{tab4.3c}. We again observe that PMLE2 has smaller bias
and standard deviation for estimating the parameters in the first
component where the mixing proportion is small, and in model III.1.1
where the two mean vectors are close. The gain is as much as 30\%
for $\sigma_{33}$.

The gains seem to disappear when the two component
mean vectors are far from each other. Nevertheless,
PMLE2 still appears to be the best estimator in terms
of both bias and standard deviation.

\vspace{3ex}%
\centerline{Tables \ref{tab4.3a}, \ref{tab4.3b} \ref{tab4.3c} about here.}
\vspace{1ex}%

We report simulation results for three models (IV.1.1, IV.2.4, IV.3.6)
in Category IV in Tables \ref{tab4.4a}, \ref{tab4.4b}, and
\ref{tab4.4c}. Again, PMLE2 has the lowest standard deviations
for estimating the parameters in the first component where the mixing
proportion is small. The comparison is the sharpest in model
IV.2.4 for $\sigma_{13}$. In contrast to the models for the other categories,
here the superiority of PMLE2 is widespread. In fact, PMLE2 is
superior for parameters in component 2, and mixed for parameters in
component 3.

We caution that even the best estimator is not necessarily a good
estimator { for trivariate mixture models}. Overall,
none of the three estimators does a great job at estimating
mixing distributions, { possibly due to their
fundamental nature, e.g., small Fisher Information
for high-dimension multivariate normal mixture models. This problem is
expected to disappear with increased sample size.}

\vspace{3ex}%
\centerline{Tables \ref{tab4.4a}, \ref{tab4.4b} \ref{tab4.4c} about here.} 
\vspace{3ex}%

{\bf Summary of the simulation results}. To conclude, {the penalized
likelihood estimators, both PMLE1 and PMLE2, are completely free
from degeneracy problems.} Moreover, 
PMLE2 has the best general performance in terms of bias
and standard deviation. This is most obvious {when the components
are not well separated}. In applications, it is unnecessary to first
judge whether it is safe to use the ratified MLE, when
a superior PMLE2 is available. Although we do not
completely dismiss the use of the ratified MLE, it is clearly
advantageous to use PMLE2 outright. We further caution against the use
of high-dimension multivariate normal mixture models in practice
when the sample size is not large. In these situations, even the
best performing estimator may not be a good estimator.

\vspace{2ex}

\newpage

{\bf Appendix}

The ordinary likelihood function is  unbounded because when the
covariance matrix of a kernel density becomes close to singular, the
likelihood contribution of the observations near its mean
vector goes to infinity. Thus, a key step in our proof is to
assess the number of such observations. In the univariate case, Chen
et al. [\ref{chen}] obtained the following result:

{\bf Lemma 1}: {\it Assume that $x_1, x_2, \ldots, x_n$ is a random
sample from a finite normal mixture distribution with density $f(x),
\, x \in R$. Let $F_n$ be the empirical distribution function and
define
$
M=\max \{ \sup_{x} f(x),  8 \},~~\mbox{and}~~ \delta_n(\sigma)= -
M \sigma \log(\sigma)+n^{-1}. ~~
$
Except for a zero-probability event not depending on $\sigma$,
we have for all large enough $n$,
\begin{enumerate}
\item[(a)]
for $\sigma$ between $\exp(-2)$ and $8/(nM)$,
$$
\sup _ \mu [F_n(\mu -\sigma \log(\sigma)) -F_n(\mu)] \leq
4 \delta_n(\sigma);
$$
\item[(b)]
 for $\sigma$ between 0 and $8/(nM)$,
$$
\sup _ \mu [F_n(\mu -\sigma \log \sigma ) -F_n(\mu)] \leq
2 n^{-1} (\log n)^2.
$$
\end{enumerate}
}

The consistency result for the multivariate normal mixture model
is built on a generalized result.
More specifically, the following lemma gives a bound
for the multivariate normal mixture model:

{\bf Lemma 2}: {\it Let $\xx_1, \xx_2, \cdots, \xx_n$ be a random
sample from a $d$-dimensional multivariate normal mixture model with
$p$ components such that its density function is given by
\[%
f(\xx,G_0)=\sum_{j=1}^{p} \pi_{j0} \varphi(\xx; ~\mu_{j0}, \Sigma_{j0}).%
\]%
Assume that all $\Sigma_{j0}$ are positive definite.
For any mean and covariance
matrix pair $(\mu, \Sigma)$ such that $|\Sigma| < \exp( - 4d)$,
except for a zero probability event not depending on $(\mu,
\Sigma)$, we have, for $n$ large enough, that
\begin{eqnarray*}
H_n {(\mu, \Sigma)}%
&=& %
\sum_{i=1}^n I\{ ({\xx_i- \mu})^\tau {\Sigma}^{-1} ({\xx_i- \mu})%
\leq%
- (\log |\Sigma|)^2 \} \\%
&\leq & %
4 (\log^2 n) I(|\Sigma| \leq \alpha_n) %
+  8 n \delta_n (|\Sigma|) I( \alpha_n \leq |\Sigma|),
\end{eqnarray*}
where
\begin{equation*}
\left\{ \begin{aligned}%
\alpha_n &= (4/Md)^{2d} n^{-2d}, \\ %
\delta_n(|\Sigma|) &=  - M |\Sigma|^{{1/2d}} \log |\Sigma| + n^{-1},
\end{aligned}%
\right .
\end{equation*}
and $M=\max\{8, \lambda_0^{-1/2} \}$ with $\lambda_0$ being the
smallest eigenvalue among those of $~\Sigma_{j0}, ~( j = 1, 2,
\ldots, p)$. %
}

{\bf Proof of Lemma 2:} Let
$~ 0 < \lambda_1 \leq \lambda_2 \leq \cdots \leq \lambda_d$ %
and $(\aa_1, \ldots , \aa_d)$ be the eigenvalues and corresponding
eigenvectors of unit length of $\Sigma$. We have that %

\begin{eqnarray*}%
&& \hspace*{-3em}
\{\xx : ({ \xx- \mu})^\tau { \Sigma}^{-1} ({ \xx- \mu})
\leq  - (\log |{\Sigma}|)^2 \}  \\
&=&%
\{%
\xx:  \sum_{j=1}^d \lambda_j^{-1} |a_j^\tau (\xx- \mu)|^2
\leq  - (\log | \Sigma |)^2 \} \\
&\subseteq&%
\{%
\xx :  |\aa_j^\tau (\xx- \mu)|  \leq  - \sqrt{\lambda_j} \log |\Sigma|, ~ j=1,\ldots,d%
\}\\
&\subseteq&%
\{%
\xx :   |\aa_1^\tau (\xx- \mu)|  \leq  - \sqrt{\lambda_1} \log |\Sigma|%
\}.
\end{eqnarray*}

Furthermore, let
\[
 Q = \{\bb_i: i=1, 2, \ldots \}
\]
be a sequence of unit vectors so that $Q$ forms a dense subset of
unit vectors in $R^d$. Hence, for any given $\aa_1$ and any
bounded subset $\mathbf{B} \in R^d$, we can find a vector $\bb$
in $Q$ such that they are arbitrarily close
so that
\[
\{%
 \xx \in \mathbf{B} :   |\aa_1^\tau ( \xx- \mu)|  \leq  - \sqrt{\lambda_1} \log |\Sigma|%
\}%
\subseteq %
\{%
\xx \in \mathbf{B}: | \bb^\tau(\xx-\mu)| \leq - \sqrt{2 \lambda_1} \log |\Sigma|%
\}.%
\]
Based on this observation, we get
\begin{eqnarray*}
\sup_{\mu} H_n {(\mu, \Sigma)}%
&=& %
\sup_{\mu} \sum_{i=1}^n I\{ (\xx_i- \mu)^\tau { \Sigma}^{-1} (\xx_i- \mu)
\leq  - (\log | \Sigma |)^2 \} \\%
&\leq&
\sup_{\bb \in Q} \sup_{\mu}
\sum_{i=1}^n I\{ | \bb^\tau (x_i- \mu)| \leq \sqrt{2 \lambda_1}|\log|\Sigma|| \}.
\end{eqnarray*}

On the other hand, given any non-random unit vector $\bb$, $\bb^\tau
\xx_i, i=1, 2, \ldots, n$ is a random sample from the univariate
normal mixture model with density
\[%
f^b(\xx)
=
 \sum_{j=1}^{p} \pi_{j0} \phi(\xx; \bb^\tau \mu_{j0}, {\bb}^\tau \Sigma_{j0} \bb ).%
\]%
We remark that since some pairs of $(\bb^\tau \mu_{j0}, {\bb}^\tau
\Sigma_{j0} \bb)$ can be equal, this univariate mixture
distribution can have fewer than $p$ components. This does not
affect the following derivation. Recall that $\lambda_0$ is the
smallest eigenvalue among those of $\Sigma_{j0}, ~j=1, \ldots, p$.
Then
\[%
\sup_{\bb \in Q} \sup_{\xx } f^b(\xx) %
\leq %
\sup_{\bb \in Q} \max \{ ({\bb}^\tau \Sigma_{j0} \bb)^{-1/2}, ~ j=1,\ldots, p \}%
=%
\lambda_0^{-\frac{1}{2}}.%
\]%
Applying Lemma 1 to the univariate
data $\bb^\tau \xx_i, i=1, \ldots, n$,
except for a zero-event not depending on
$\Sigma$, as $n \to \infty$, we have
\begin{eqnarray*}
&&
\sup_{\mu} \sum_{i=1}^n I\{ | \bb^\tau (\xx_i - \mu)|
\leq
\sqrt{\lambda_1}|\log |\Sigma|| \} \\%
&\leq& %
4 (\log^2 n) I(|\Sigma| \leq \alpha_n) %
+  8 n \delta_n (|\Sigma|) I( \alpha_n \leq |\Sigma|).%
\end{eqnarray*}
The conclusion of the lemma simply claims that
the above inequality is true over all $\bb \in Q$
with only a zero-probability-event exception.
The zero-probability claim remains true because
$Q$ is countable.

\vspace{2ex}%
\noindent %
{\bf Proof of Theorem 1}: We give a proof for the case $p=2$;
the proof for the general case is similar. Let $\Gamma$ be the
parameter space for $G$ and define
\begin{eqnarray*}
\Gamma_1 &=& \{G\in \Gamma \, : \, |\Sigma_1| \leq |\Sigma_2| \leq \varepsilon_0\} \\ %
\Gamma_2 &=& \{G\in \Gamma \, : \, |\Sigma_1| \leq\tau_0, |\Sigma_2| \geq\varepsilon_0\} \\ %
\Gamma_3 &=& \Gamma-(\Gamma_1 \cup \Gamma_2) %
\end{eqnarray*}
where $\varepsilon_0>\tau_0>0$ are two small positive
constants to be specified soon. The first subspace represents
the case where the two components have nearly singular covariance
matrices. Hence the observations inside the small ellipse centered
at the mean parameter make a large contribution to the log
likelihood function.

Let $K_0 = E \{ \log f(X; G_0) \}$. 
The constants
$\varepsilon_0$, $\tau_0$
must satisfy the following four conditions:\\ %
\begin{enumerate}
\item[1:] %
$0<\varepsilon_0<\exp\{ - 4d \}$;
\item[2:]%
$- \log \varepsilon_0-(\log \varepsilon_0)^2 \leq  4(K_0-2)$;
\item[3:]%
$16 M \varepsilon_0^{1/2d} (\log \varepsilon_0)^2 \leq 1 $;
\item[4:]%
$16 M d \tau_0(\log\tau_0)^2 \leq \frac{2}{5}\delta_0$;
\end{enumerate}
for some $\delta_0 > 0$ to be specified. The existence
of $\varepsilon_0, \, \tau_0$ is obvious.

We proceed with the proof in three steps.

\noindent {\bf Step 1.}
For any $G \in \Gamma_1$, we show that almost surely,
$$
\sup_{\Gamma_1}p l_n(G)-p l_n(G_0) \to  -\infty.
$$
Define two index sets
\begin{eqnarray*}
 A
 &=&
  \{i:(x_i-\mu_1)^\tau \Sigma_1^{-1}(x_i-\mu_1)\leq(\log|\Sigma_1|)^2\},\\
 B
 &=&
  \{i:(x_i-\mu_2)^\tau \Sigma_2^{-1}(x_i-\mu_2)\leq(\log|\Sigma_2|)^2\},
\end{eqnarray*}
and for any index set $S \in \{1,2,\ldots,n\}$, denote
$$
l_n(G; ~ S) = \sum_{i \in S} \log f(X_i,~G).
$$
We can write %
$l_n(G)=l_n(G;~ A)+l_n(G;~ A^c B)+l_n(G;~ A^c B^c)$, where $A^c$
and $ B^c$ are the complement sets of $A$ and $B$ respectively. For
any index set $S$, denote $n(S)$ as its cardinality. It is easy to
see that
$$
l_n(G;~ A) \leq  n(A) \log|\Sigma_1|^{-\frac{1}{2}},
$$
$$
l_n(G;~ B) \leq n(B) \log|\Sigma_2|^{-\frac{1}{2}}.
$$
Applying Lemma 2 
to $n(A)$ and $n(B)$, noting that
$|\Sigma_1 | \leq \epsilon_0$ for $G$ in $\Gamma_1$, and C3 on the
penalty function, we find that
$$
l_n(G;~ A)+ \widetilde{p}_n(\Sigma_1) %
\leq %
16 d \log n %
+ %
8 M \varepsilon_0^{\frac{1}{2d}} (\log \varepsilon_0)^2 n
$$
$$
l_n(G;~ A^c B)+\widetilde{p}_n(\Sigma_2)  %
\leq %
16 d \log n %
+ %
8 M  \varepsilon_0^{\frac{1}{2d}} (\log \varepsilon_0)^2 n.
$$
The key point underlying the above two inequalities is that they
are bounded by an arbitrarily small fraction of $n$.
Further, for observations away from $\mu_1$ and $\mu_2$,
we have
\begin{eqnarray*}
&&l_n(G;~ A^c B^c) \\%
&\leq& %
\sum_{i \in A^c B^c}
\log [ \pi_1 \exp\{\log|\Sigma_1|^{-\frac{1}{2}}
-
\frac{1}{2}(\log |\Sigma_1|)^2\}%
+%
\pi_2 \exp\{\log|\Sigma_2|^{-\frac{1}{2}}
-
\frac{1}{2}(\log|\Sigma_2|)^2\} ]
\\
&\leq& %
\sum_{i\in A^c B^c}
\{-\frac{1}{2}\log \varepsilon_0 - \frac{1}{2}(\log \varepsilon_0)^2\} \\%
&\leq& %
n(K_0-2)%
\end{eqnarray*}
The last line in the above derivation is obtained
by choosing a small enough
$\epsilon_0$ as specified earlier.
Combining these inequalities,
we get $p l_n(G)\leq n(K_0-1)$, and hence almost surely
\[
\sup_{\Gamma_1}pl_n(G)-pl_n(G_0)\leq -n + 16 d \log n.
\] %
That is,
$$
\sup_{\Gamma_1}p l_n(G)-p l_n(G_0)\rightarrow -\infty
$$
almost surely which completes the first step.

\noindent {\bf Step 2.}
For $G \in \Gamma_2,$ we also show that almost surely
$$
 \sup_{\Gamma_2}p l_n(G)-p l_n(G_0)\rightarrow -\infty.
$$
Recall that for each $i \in A$,
$(\xx_i - \mu_1)^\tau \Sigma_1^{-1} (\xx_i - \mu_1)$
is bounded by $(\log \Sigma_1)^2$.
Hence, it is easy to verify that for $i \in A$,
$$%
\varphi({\xx_i}; \mu_1,\Sigma_1) %
\leq%
|\Sigma_1|^{-1/2}
\exp\{-\frac{1}{4}({\xx_i}-\mu_1)^\tau \Sigma_1^{-1} (\xx_i-\mu_1)\}.
$$
For $i \not \in A$,
$$%
\varphi(\xx_i; \mu_1,\Sigma_1)
\leq
\exp\{-\frac{1}{4}({\xx_i}-\mu_1)^T \Sigma_1^{-1} ({\xx_i}-\mu_1)\}. %
$$
Therefore, letting (not a density itself)
$$%
g(\xx; G)%
=%
\pi_1 \exp\{-\frac{1}{4}(\xx-\mu_1)^T\Sigma_{1}^{-1}(\xx-\mu_1)\}%
+%
\pi_2 \varphi(\xx;\mu_2,\Sigma_2),%
$$
we have
\[
\log f(\xx_i; G) \leq \log g(\xx_i; G) +I(i \in A) \log |\Sigma_1|^{-1/2}.
\]
Hence, we get
\[%
l_n(G;~ A) \leq %
n(A)\log|\Sigma_1|^{-\frac{1}{2}} %
+ %
\sum_{i=1}^n g({\xx_i}; G).
\]%

It is obvious that for any $G \in \Gamma_2$, (a)
$E_0 \left\{ \log g(X; G) / f(X; G_0) \right\}< 0$ by Jensen's inequality
and the fact that the integration of $g(\xx, G)$ is less than 1;
(b) $g(\xx; G) \leq \varepsilon_0^{-1}$ by the definition of $\Gamma_2$.
Hence for each given $G \in \Gamma_2$, by the law of large numbers,
\begin{equation*}
\label{eqn31}
\frac{1}{n} \sum_{i=1}^n
\log \{ g(X_i; G)/f(X_i; G_0)\}
\to
E \{ g(X; G)/f(X; G_0) \} < 0.
\end{equation*}
For each fixed $\xx$, we can extend the definition of $g(\xx; G)$ in
$G$ onto the compacted $\Gamma_2$ while maintaining properties (a)
and (b) and its continuity in $G$. Thus, a classical technique as
in [\ref{wald}] can be readily employed to show that as $n \to
\infty$,
\begin{equation}\label{eqn new31}
\sup_{G \in \Gamma_2}
\left\{ \frac{1}{n} \sum_{i=1}^n \log \left(\frac{g(X_i;G)}{f(X_i;G_0)} \right) \right\}
\to - \delta(\tau_0)<0
\end{equation}
for some decreasing function $\delta(\tau_0)$. Hence,
it is possible to choose a small enough $\tau_0 \leq \epsilon_0$,
such that
$$%
\aligned
& \sup_{\Gamma_2} pl_n(G) - pl_n(G_0) \\%
& \leq %
\sup_{\Gamma_2}\{n(A)\log|\Sigma_1|^{-\frac{1}{2}} + p_n(G) \}
+%
\sup_{{\Gamma}_2} \sum_{i=1}^n \log \left\{ \frac{g(X_i,G)} {f(X_i,G_0)} \right\} \\
& \leq %
8 M \tau_0 (\log \tau_0)^2 n - \frac{9}{10}\delta(\epsilon_0) n\\
& \leq%
-\frac{1}{2} \delta(\epsilon_0) n.
\endaligned%
$$
The first term of the third line above is from the assessment of $n(A)$,
C3 on $p_n(G)$. Note also that $p_n(G_0) = o(n)$.
Therefore, almost surely,
$$%
\sup_{\Gamma_2}p l_n(G)-p l_n(G_0)\rightarrow -\infty.
$$

\noindent {\bf Step 3.} From the above two steps, we know that
$\tilde G_n \in \Gamma_3$ with probability 1. At the same time,
when $G \in \Gamma_3$, we have $p_n(G) = o(1)$.
By the definition of the maximum penalized likelihood estimator,
we have
\begin{equation}
\label{last} l_n (\tilde G_n) -
l_n(G_0) \geq p_n(G) - p_n(G_0) = o(1).%
\end{equation}

Since the parameter space $\Gamma_3$ is now completely
regular, an estimator with property (\ref{last}) is easily shown
to be consistent by the classical technique [\ref{wald}]
even with a penalty of size $o(n)$.
\hfill{$\Box$}

\vspace{2ex} %
\noindent %
{\bf Proof of Theorem 3}: When $p_0 < p < \infty$, we cannot expect
that every part of $G$ converges to that of $G_0$. Instead, we
measure their difference as two distributions. Let
\begin{equation*}
H(G, G_0) = \int_{{\cal R}^d \times {\mathcal{A}}}
|G(\mathbf{\lambda})-G_0(\mathbf{\lambda})|
\exp\{-|\mathbf{\lambda}| \} d \mathbf{\lambda}
\end{equation*}
where
\[
\mathbf{\lambda} = (\mu_1, \mu_2,..., \mu_d, \sigma_{11},
\sigma_{12}, \sigma_{22},..., \sigma_{dd}) \in {\cal R}^d \times
{\mathcal{A}},
\]
\[
|\mathbf{\lambda}| = \sum_{j=1}^d |\mu_j|  + \sum_{i=1}^d
\sum_{j=1}^i |\sigma_{ij}|,
\]
and ${\mathcal{A}}$ is a subset of ${\cal R}^{d \times (d+1)/2}$
containing all eligible combinations of ${d \times (d+1)/2}$ real
numbers which  form a symmetric positive definite matrix. It is
well known that ${\mathcal{A}}$ is an open connected subset of
${\cal R}^{d \times (d+1)/2}$ and is regular enough although it
may not be easy to visualize its shape.

It can be shown that $ H(G_n, G_0) \to 0 $ implies $ G_n \to G_0 $
in distribution. An estimator $\tilde G_n$ is strongly consistent if
$H(\tilde G_n, G_0) \to 0$ almost surely.

Again, for the sake of clarity, we consider only the special case
with $p=2, p_0=1$, that is, to fit a non-mixture multivariate
normal model with a two-component multivariate normal mixture
 model. The extension of our proof to general situations is
 straightforward and the major hurdle is merely a
complicated presentation.
Most intermediate conclusions in the proof of consistency of the
PMLE when $p=p_0 = 2$ are still applicable; some need minor
changes. We use many of these results and notations to establish a
brief proof.

For an arbitrarily small positive number $\delta$, define \( {\cal
H}(\delta) = \{G: G \in \Gamma,  H(G, G_0) \geq \delta \}. \) That
is, ${\cal H}(\delta)$ contains all mixing distributions with up
to $p$ components that are at least $\delta >0$ distance from the
true mixing distribution $G_0$.

Since $G_0 \not \in {\cal H}(\delta)$, we have \( E [\log \{g(X;
G)/f(X; G_0)\}] < 0 \) for any $G \in {\cal H}(\delta)$. 
Thus, (\ref{eqn new31}) remains valid after being slightly revised
as follows:
\[
\sup_{G \in {\cal H}(\delta)\cap \Gamma_2} n^{-1} \sum_{i=1}^n
\log \{g(X_i; G)/f(X_i; G_0)\} \to - \eta(\tau)
\]
for some positive $\eta(\tau)$ depending on $\Gamma_2$. Because of
this, the derivations in the proof of Theorem 1 still apply after
$\Gamma_k$ is replaced by ${\cal H}(\delta) \cap \Gamma_k$
($k=1,2$). That is, with proper choice of $\epsilon_{0}$ and
$\tau_0$, we similarly get \( \sup_{G \in {\cal H}(\delta)\cap
\Gamma_k} pl_n(G) - pl_n(G_0) \to -\infty \) for $k=1, 2$.

With what we have proved, it is seen that the penalized maximum
likelihood estimator of $G$, $\tilde G_n$, must almost surely belong
to ${\cal H}^c(\delta)\cup \Gamma_{3}$, where ${\cal H}^c(\delta)$
is the complement of ${\cal H}(\delta)$. Since $\delta$ is
arbitrarily small, $\tilde G_n \in {\cal H}^c(\delta)$ implies
$H(\tilde G_n, G_0) \to 0$. On the other hand, $\tilde G_n \in
\Gamma_{3}$ is equivalent to putting a positive lower bound on the
component variances, which also implies $H(\tilde G_n, G_0) \to 0$
by [\ref{kie}].
That is, consistency of the PMLE is
also true when $p=2$ but $p_0=1$.

A generalization of the above derivation leads to the conclusion
of Theorem 3.

\newpage

\begin{table}[p]
\caption{Number of Degeneracies} 
\begin{center} 
\vspace{1ex}
\renewcommand{\arraystretch}{1.0}
\begin{tabular}{|c|r|r|r|r|r|r|}
\hline
Mean.Var.Config &    1 &    2&     3&      4 &    5 &  6 \\ \hline%
 & \multicolumn{6}{c|}{2-component bivariate normal mixture} \\
 \hline
near &  0      &  11  &   19  &   5     & 40   &  8    \\%
moderate &  1911   &  3256&   441 &   6     & 2523 &  157  \\%
distant &  4997   &  4998&   4966&   4782  & 4998 &  4943  \\%
\hline
 & \multicolumn{6}{c|}{3-component bivariate normal mixture}\\
  \hline
straight &  3049   &  5058&   4947&   1998  & 2306 &  2491 \\%
acute &  2888   &  4505&   4812&   4052  & 4057 &  4561 \\%
obtuse &  3253   &  4980&   4983&   2885  & 3022 &  3511  \\%
\hline
 & \multicolumn{6}{c|}{2-component trivariate normal mixture}\\
  \hline
near &  1      &  4872&   5003&   4866  & 4961 &  1466 \\%
moderate &  4011   &  5000&   5001&   5000  & 5000 &  4900 \\%
distant &  5000   &  5000&   5000&   5000  & 5000 &  5000  \\%
\hline
 & \multicolumn{6}{c|}{3-component trivariate normal mixture}\\
  \hline
straight &  5009   &  5010&   5002&   5002  & 5000 &  5000 \\%
acute &  5006   &  5034&   5000&   5002  & 5000 &  5000 \\%
obtuse &  5009   &  5038&   5002&   5004  & 5000 &  5001  \\%
\hline %
\end{tabular}\\
\end{center}
\end{table}

\begin{table}[p]
\caption{Bias (std) under 2-component bivariate normal mixture
models.}
\label{tab4.1} %
\begin{center}
\vspace{1ex}
\renewcommand{\arraystretch}{1.0}
\begin{tabular}{|l|c|c|c|}
  \hline
& {MLE} & {PMLE1} & {PMLE2} \\ \hline%
&\multicolumn{3}{c|}{Model I.1.1, component 1} \\ \hline%
 $\pi = 0.3    $& -0.03 (0.11) & -0.02 (0.11) & -0.01 (0.10)\\
 $\mu_1= 0     $& -0.16 (0.53) & -0.16 (0.53) & -0.13 (0.50)\\
 $\mu_2= -1    $& ~0.72 (1.17) & ~0.72 (1.17) & ~0.71 (1.14)\\
 $\sigma_{11}=1$& -0.14 (0.41) & -0.14 (0.40) & -0.13 (0.37)\\ %
 $\sigma_{12}=0$& -0.01 (0.39) & ~0.00 (0.38) & ~0.00 (0.34)\\ %
 $\sigma_{22}=1$& -0.03 (0.71) & -0.03 (0.70) & -0.01 (0.64)\\
\hline
&\multicolumn{3}{c|}{Model I.1.1, component 2} \\ \hline%
 $\pi_2 = 0.7  $& ~0.03 (0.11) & ~0.02 (0.11) & ~0.01 (0.10)\\
 $\mu_1=0      $& ~0.04 (0.19) & ~0.04 (0.19) & ~0.04 (0.19)\\
 $\mu_2=1      $& -0.39 (0.47) & -0.39 (0.47) & -0.37 (0.48)\\
 $\sigma_{11}=1$& -0.07 (0.18) & -0.07 (0.18) & -0.07 (0.18)\\
 $\sigma_{12}=0$& ~0.00 (0.19) & ~0.00 (0.19) & ~0.00 (0.19)\\
 $\sigma_{22}=1$& ~0.33 (0.44) & ~0.33 (0.44) & ~0.30 (0.43)\\
\hline
&\multicolumn{3}{c|}{Model I.2.4, component 1} \\ \hline%
 $\pi_1 = 0.3  $& 0.00 (0.03)&  0.00 (0.03)&  0.00 (0.03)\\
 $\mu_1 = 0    $&-0.02 (0.28)& -0.02 (0.28)& -0.02 (0.28)\\
 $\mu_2 = -3   $&-0.01 (0.13)& -0.01 (0.13)& -0.01 (0.13)\\
 $\sigma_{11}=5$&-0.04 (0.93)& -0.04 (0.93)& -0.04 (0.93)\\
 $\sigma_{12}=0$& 0.00 (0.30)&  0.00 (0.30)&  0.00 (0.30)\\
 $\sigma_{22}=1$&-0.02 (0.19)& -0.02 (0.19)&  0.00 (0.19)\\
\hline
&\multicolumn{3}{c|}{Model I.2.4, component 2} \\ \hline%
 $\pi_2 = 0.7  $& 0.00 (0.03)&  0.00 (0.03)&  0.00 (0.03)\\
 $\mu_1=0      $& 0.00 (0.09)&  0.00 (0.09)&  0.00 (0.09)\\
 $\mu_2=3      $& 0.00 (0.09)&  0.00 (0.09)&  0.00 (0.09)\\
 $\sigma_{11}=1$&-0.01 (0.12)& -0.01 (0.12)& -0.01 (0.12)\\
 $\sigma_{12}=0$& 0.00 (0.08)&  0.00 (0.08)&  0.00 (0.08)\\
 $\sigma_{22}=1$& 0.00 (0.12)&  0.00 (0.12)&  0.00 (0.12)\\
 \hline
\end{tabular}
\end{center}
\end{table}

\begin{table}[p]
\caption{Bias (std) under 3-component bivariate
normal mixture models.}
\label{tab4.2a}
\begin{center}
\vspace{1ex}
\renewcommand{\arraystretch}{1.0}
\begin{tabular}{|l|c|c|c|}
  \hline
& {MLE} & {PMLE1} & {PMLE2} \\ \hline%
&\multicolumn{3}{c|}{Model II.1.1, component 1} \\ \hline%
 $\pi = 0.15$     & -0.10 (0.06) & -0.08 (0.07)& -0.04 (0.07)\\
 $\mu_1= 0$       & ~0.69 (1.15) & ~0.58 (1.28)& ~0.25 (1.01)\\
 $\mu_2= -2$      & ~1.17 (2.48) & ~1.15 (2.32)& ~1.24 (1.94)\\
 $\sigma_{11} = 1$& -0.33 (0.91) & -0.46 (0.60)&-0.33 (0.52)\\ %
 $\sigma_{12} = 0$& -0.04 (0.54) & -0.02 (0.46)&~0.02 (0.48)\\ %
 $\sigma_{22} = 1$& -0.22 (1.16) & -0.22 (1.01)&~0.12 (1.01)\\
 \hline
 &\multicolumn{3}{c|}{Model II.1.1, component 2} \\ \hline%
 $\pi_2 = 0.35$    & -0.02 (0.10) & -0.02 (0.10) & -0.03 (0.08)\\
 $\mu_1 = 0$       & -0.10 (0.39) & -0.08 (0.38) & -0.06 (0.39)\\
 $\mu_2 = 0$       & ~0.61 (1.54) & ~0.63 (1.53) & ~0.56 (1.44)\\
 $\sigma_{11} = 1$ & -0.13 (0.29) & -0.13 (0.30) & -0.14 (0.31)\\
 $\sigma_{12} = 0$ & ~0.02 (0.32) & ~0.01 (0.33) & ~0.02 (0.34)\\
 $\sigma_{22} = 1$ & ~0.24 (0.70) & ~0.20 (0.71) & ~0.22 (0.69)\\
 \hline
 &\multicolumn{3}{c|}{Model II.1.1, component 3} \\ \hline%
 $\pi_3 = 0.5$     & ~0.11 (0.11) & ~0.10 (0.12) & ~0.06 (0.10)\\
 $\mu_1 = 0$       & ~0.02 (0.20) & ~0.01 (0.21) & ~0.01 (0.24)\\
 $\mu_2 = 2$       & -1.23 (0.90) & -1.16 (0.89) & -1.02 (0.89)\\
 $\sigma_{11} = 1$ & -0.08 (0.16) & -0.08 (0.17) & -0.10 (0.19)\\
 $\sigma_{12} = 0$ & ~0.03 (0.26) & ~0.03 (0.27) & ~0.00 (0.28)\\
 $\sigma_{22} = 1$ & ~0.86 (0.68) & ~0.81 (0.70) & ~0.65 (0.67)\\
 \hline
\end{tabular}
\end{center}
\end{table}

\begin{table}[p]
\caption{Bias (std) under 3-component bivariate normal mixture
models.} \label{tab4.2b}
\begin{center}
\vspace{1ex}
\renewcommand{\arraystretch}{1.0}
\begin{tabular}{|l|c|c|c|}
  \hline
 &\multicolumn{3}{c|}{Model II.2.4, component 1} \\ \hline%
 $\pi_1 = 0.15 $ & 0.00 (0.04)& 0.01 (0.04)& 0.01 (0.03)\\
 $\mu_1 = 0 $    & 0.23 (0.86)& 0.18 (0.74)& 0.19 (0.72)\\
 $\mu_2 = -2$    & 0.12 (0.83)& 0.11 (0.63)& 0.11 (0.54)\\
 $\sigma_{11}=1$ & 0.07 (0.69)& 0.06 (0.60)& 0.10 (0.59)\\
 $\sigma_{12}=0$ &-0.05 (0.54)&-0.03 (0.40)&-0.04 (0.38)\\
 $\sigma_{22}=1$ & 0.17 (0.99)& 0.18 (0.95)& 0.20 (0.90)\\
\hline
 &\multicolumn{3}{c|}{Model II.2.4, component 2} \\ \hline%
 $\pi_2 = 0.35 $&-0.01 (0.05)&-0.01 (0.05)&-0.01 (0.05)\\
 $\mu_1 = 3 $   &-0.43 (1.12)&-0.40 (1.09)&-0.38 (1.08)\\
 $\mu_2 = 0 $   & 0.15 (0.82)& 0.14 (0.80)& 0.13 (0.79)\\
 $\sigma_{11}=1$& 0.37 (1.12)& 0.33 (1.05)& 0.31 (1.03)\\
 $\sigma_{12}=0$&-0.01 (0.35)&-0.02 (0.34)&-0.03 (0.37)\\
 $\sigma_{22}=5$&-0.69 (1.60)&-0.65 (1.57)&-0.62 (1.55)\\
\hline               %
 &\multicolumn{3}{c|}{Model II.2.4, component 3} \\ \hline%
 $\pi_3 = 0.5  $& 0.00 (0.05)& 0.00 (0.05)& 0.00 (0.05)\\
 $\mu_1 = 0 $   & 0.33 (0.88)& 0.31 (0.88)& 0.30 (0.87)\\
 $\mu_2 = 2 $   &-0.19 (0.57)&-0.17 (0.53)&-0.16 (0.51)\\
 $\sigma_{11}=5$&-0.38 (1.31)&-0.36 (1.31)&-0.36 (1.30)\\
 $\sigma_{12}=0$& 0.00 (0.28)&-0.01 (0.26)&-0.01 (0.27)\\
 $\sigma_{22}=1$& 0.37 (1.15)& 0.34 (1.11)& 0.33 (1.08)\\
\hline               %
\end{tabular}
\end{center}
\end{table}

\begin{table}[p]
\caption{Bias (std) under 2-component trivariate
normal mixture models.}
\label{tab4.3a}
\begin{center}
\vspace{1ex}
\renewcommand{\arraystretch}{1.0}
\begin{tabular}{|l|c|c|c|}
  \hline
& {MLE} & {PMLE1} & {PMLE2} \\ \hline%
&\multicolumn{3}{c|}{Model III.1.1, component 1} \\ \hline%
$\pi_1 = 0.3 $    & -0.09 (0.15)& -0.08 (0.15)& -0.05 (0.14)\\
$\mu_1 = 0$       & -0.28 (0.61)& -0.26 (0.58)& -0.17 (0.51)\\
$\mu_2 = 0$       & -0.15 (0.58)& -0.14 (0.57)& -0.09 (0.52)\\
$\mu_3 =-1$       & ~0.52 (0.09)& ~0.54 (0.11)& ~0.61 (0.09)\\
$\sigma_{11} = 1$ & -0.12 (0.47)& -0.11 (0.46)& -0.11 (0.36)\\
$\sigma_{12} = 0$ & -0.01 (0.38)& ~0.00 (0.35)& ~0.02 (0.27)\\
$\sigma_{13} = 0$ & -0.10 (0.48)& -0.10 (0.47)& -0.07 (0.37)\\
$\sigma_{22} = 1$ & -0.09 (0.56)& -0.11 (0.47)& -0.13 (0.36)\\
$\sigma_{23} = 0$ & -0.04 (0.49)& -0.02 (0.47)& -0.01 (0.37)\\
$\sigma_{33} = 1$ & ~0.22 (0.91)& ~0.18 (0.83)& ~0.12 (0.66)\\
\hline              %
&\multicolumn{3}{c|}{Model III.1.1, component 2} \\ \hline%
$ \pi_2 = 0.7$    & ~0.09 (0.15)& ~0.08 (0.15)& ~0.05 (0.14)\\
$\mu_1 = 0$       & ~0.01 (0.15)& ~0.01 (0.15)& ~0.01 (0.16)\\
$\mu_2 = 0$       & ~0.02 (0.15)& ~0.02 (0.15)& ~0.02 (0.17)\\
$\mu_3 = 1$       & -0.45 (0.41)& -0.44 (0.41)& -0.42 (0.44)\\
$\sigma_{11} = 1$ & -0.05 (0.13)& -0.05 (0.13)& -0.05 (0.14)\\
$\sigma_{12} = 0$ & ~0.00 (0.10)& ~0.00 (0.10)& ~0.00 (0.10)\\
$\sigma_{13} = 0$ & -0.02 (0.13)& -0.02 (0.13)& -0.02 (0.14)\\
$\sigma_{22} = 1$ & ~0.03 (0.13)& -0.03 (0.13)& -0.04 (0.14)\\
$\sigma_{23} = 0$ & ~0.01 (0.14)& ~0.01 (0.14)& ~0.01 (0.15)\\
$\sigma_{33} = 1$ & ~0.44 (0.38)& ~0.43 (0.38)& ~0.39 (0.39)\\
\hline
\end{tabular}
\end{center}
\end{table}

\begin{table}[p]
\caption{Bias (std) under 2-component trivariate
normal mixture models.}
\label{tab4.3b}
\begin{center}
\vspace{1ex}
\renewcommand{\arraystretch}{1.0}
\begin{tabular}{|l|c|c|c|}
  \hline
&\multicolumn{3}{c|}{Model III.2.4, component 1} \\ \hline%
$\pi_1 = 0.3$     & 0.00  (0.04) & 0.00  (0.04) & 0.00  (0.04)\\
$\mu_1 = 0$       & 0.01  (0.13) & 0.01  (0.13) & 0.01  (0.13)\\
$\mu_2 = 0$       & 0.01  (0.22) & 0.01  (0.22) & 0.01  (0.22)\\
$\mu_3 =-3$       & -0.03 (0.52) & -0.03 (0.52) & -0.04 (0.52)\\
$\sigma_{11} = 1$ & -0.01 (0.17) & -0.01 (0.17) & -0.01 (0.17)\\
$\sigma_{12} = 0$ & -0.01 (0.20) & -0.01 (0.20) & -0.01 (0.19)\\
$\sigma_{13} = 0$ & 0.03  (0.45) & 0.03  (0.45) & 0.03  (0.45)\\
$\sigma_{22} = 3$ & -0.05 (0.49) & -0.05 (0.49) & -0.04 (0.49)\\
$\sigma_{23} = 0$ & 0.00  (0.75) & 0.00  (0.75) & 0.01  (0.75)\\
$\sigma_{33} =10$ & -0.36 (2.10) & -0.36 (2.11) & -0.38 (2.09)\\
\hline           %
&\multicolumn{3}{c|}{Model III.2.4, component 2} \\ \hline%
$\pi_2 = 0.7 $        & 0.00  (0.04) & 0.00  (0.04) & 0.00  (0.04)\\
$\mu_1 = 0$           & 0.00  (0.15) & 0.00  (0.15) & 0.00  (0.15)\\
$\mu_2 = 0$           & -0.01 (0.19) & -0.01 (0.19) & -0.01 (0.19)\\
$\mu_3 = 3$           & -0.01 (0.11) & -0.01 (0.11) & -0.01 (0.11)\\
$\sigma_{11}=4.87 $   & -0.03 (0.47) & -0.03 (0.48) & -0.03 (0.47)\\
$\sigma_{12}=-3.23$   & 0.03  (0.49) & 0.03  (0.49) & 0.03  (0.48)\\
$\sigma_{13}=-0.5 $   & 0.01  (0.23) & 0.01  (0.23) & 0.01  (0.23)\\
$\sigma_{22}=7.2  $   & -0.07 (0.71) & -0.07 (0.72) & -0.07 (0.71)\\
$\sigma_{23}=2.16 $   & -0.02 (0.30) & -0.02 (0.30) & -0.02 (0.30)\\
$\sigma_{33}=1.94 $   & -0.01 (0.22) & -0.01 (0.22) & 0.00  (0.22)\\
\hline
\end{tabular}
\end{center}
\end{table}

\begin{table}[p]
\caption{Bias (std) under 2-component trivariate
normal mixture models.}
\label{tab4.3c}
\begin{center}
\vspace{1ex}
\renewcommand{\arraystretch}{1.0}
\begin{tabular}{|l|c|c|c|}
  \hline
&\multicolumn{3}{c|}{Model III.3.6, component 1} \\ \hline%
$\pi_1  = 0.3  $  & 0.00  (0.03) & 0.00  (0.03) & 0.00  (0.03)\\
$\mu_1=0       $  & 0.00  (0.10) & 0.00  (0.10) & 0.00  (0.10)\\
$\mu_2=0       $  & 0.01  (0.19) & 0.01  (0.19) & 0.00  (0.19)\\
$\mu_3=-5      $  & 0.01  (0.37) & 0.01  (0.37) & 0.01  (0.37)\\
$\sigma_{11}= 1$  & -0.01 (0.15) & -0.01 (0.15) & -0.01 (0.15)\\
$\sigma_{12}= 0$  & 0.01  (0.18) & 0.01  (0.18) & 0.01  (0.18)\\
$\sigma_{13}= 0$  & 0.02  (0.36) & 0.02  (0.36) & 0.02  (0.36)\\
$\sigma_{22}= 3$  & -0.05 (0.45) & -0.05 (0.45) & -0.04 (0.45)\\
$\sigma_{23}= 0$  & -0.02 (0.64) & -0.02 (0.64) & -0.02 (0.64)\\
$\sigma_{33}=10$  & -0.06 (1.81) & -0.06 (1.81) & -0.06 (1.80)\\
\hline           %
&\multicolumn{3}{c|}{Model III.3.6, component 2} \\ \hline%
$\pi_2     =0.7$      & 0.00  (0.03) & 0.00  (0.03) & 0.00  (0.03)\\
$\mu_1      =0 $      & 0.00  (0.15) & 0.00  (0.15) & 0.00  (0.15)\\
$\mu_2      =0 $      & 0.00  (0.19) & 0.00  (0.19) & 0.00  (0.19)\\
$\mu_3    =  5 $      & 0.00  (0.10) & 0.00  (0.10) & 0.00  (0.10)\\
$\sigma_{11} = 4.87$  & -0.05 (0.46) & -0.05 (0.46) & -0.05 (0.46)\\
$\sigma_{12} =3.23  $ & -0.03 (0.46) & -0.03 (0.46) & -0.03 (0.46)\\
$\sigma_{13} =-0.5  $ & 0.00  (0.22) & 0.00  (0.22) & 0.00  (0.22)\\
$\sigma_{22} =7.2   $ & -0.02 (0.70) & -0.02 (0.70) & -0.03 (0.70)\\
$\sigma_{23} =-2.16 $ & -0.01 (0.29) & -0.01 (0.29) & -0.01 (0.29)\\
$\sigma_{33}=1.94  $  & -0.01 (0.20) & -0.01 (0.20) & 0.00  (0.20)\\
\hline
\end{tabular}
\end{center}
\end{table}

\begin{table}[p]
\caption{Bias (std) under 3-component trivariate
normal mixture models.}
\label{tab4.4a}
\begin{center}
\vspace{1ex}
\renewcommand{\arraystretch}{1.0}
\begin{tabular}{|l|c|c|c|}
  \hline
& {MLE} & {PMLE1} & {PMLE2} \\ \hline%
&\multicolumn{3}{c|}{Model IV.1.1, component 1} \\ \hline%
$\pi_1 =0.15   $& -0.05 (0.07)&-0.06 (0.07)&-0.01 (0.07)\\
$\mu_1    =0$   &  0.10 (0.64)& 0.28 (0.97)& 0.12 (0.69)\\
$\mu_2    =0$   & -0.08 (0.64)& 0.11 (0.97)&-0.04 (0.65)\\
$\mu_3   =-2$   &  3.07 (2.16)& 2.65 (2.17)& 2.16 (1.89)\\
$\sigma_{11} =1$& -0.05 (0.73)&-0.25 (0.63)&-0.19 (0.47)\\
$\sigma_{12} =0$&  0.07 (0.50)& 0.05 (0.40)& 0.04 (0.35)\\
$\sigma_{13} =0$& -0.01 (0.58)& 0.00 (0.51)& 0.00 (0.48)\\
$\sigma_{22} =1$& -0.04 (0.74)&-0.23 (0.63)&-0.16 (0.47)\\
$\sigma_{23} =0$&  0.03 (0.51)& 0.03 (0.47)& 0.04 (0.43)\\
$\sigma_{33} =1$& -0.01 (1.16)& 0.01 (1.19)& 0.31 (1.05)\\
\hline              %
&\multicolumn{3}{c|}{Model IV.1.1, component 2} \\ \hline%
$\pi_2 =0.35   $& -0.05 (0.09)&-0.07 (0.11)&-0.05 (0.09)\\
$\mu_1    =0   $& -0.05 (0.33)&-0.10 (0.45)&-0.02 (0.37)\\
$\mu_2    =0$   &  0.04 (0.33)&-0.02 (0.43)& 0.01 (0.34)\\
$\mu_3    =0$   &  0.00 (1.47)& 0.02 (1.52)& 0.26 (1.42)\\
$\sigma_{11} =1$& -0.09 (0.26)&-0.12 (0.32)&-0.11 (0.29)\\
$\sigma_{12} =0$&  0.02 (0.20)& 0.01 (0.23)& 0.02 (0.21)\\
$\sigma_{13} =0$& -0.05 (0.32)&-0.05 (0.41)&-0.03 (0.35)\\
$\sigma_{22} =1$& -0.09 (0.28)&-0.11 (0.30)&-0.11 (0.28)\\
$\sigma_{23} =0$&  0.02 (0.33)&-0.01 (0.37)& 0.01 (0.33)\\
$\sigma_{33} =1$&  0.46 (0.83)& 0.48 (0.93)& 0.46 (0.84)\\
\hline              %
&\multicolumn{3}{c|}{Model IV.1.1, component 3} \\ \hline%
$\pi_3 =0.5    $&  0.10 (0.12)& 0.13 (0.15)& 0.06 (0.12)\\
$\mu_1    =0$   &  0.01 (0.19)& 0.00 (0.18)& 0.00 (0.21)\\
$\mu_2    =0$   & -0.01 (0.18)&-0.01 (0.17)& 0.00 (0.21)\\
$\mu_3    =2$   & -0.96 (0.81)&-1.00 (0.79)&-0.97 (0.86)\\
$\sigma_{11} =1$& -0.07 (0.17)&-0.07 (0.17)&-0.08 (0.19)\\
$\sigma_{12} =0$&  0.01 (0.12)& 0.00 (0.11)& 0.01 (0.13)\\
$\sigma_{13} =0$& -0.04 (0.22)&-0.04 (0.22)&-0.04 (0.24)\\
$\sigma_{22} =1$& -0.06 (0.16)&-0.06 (0.16)&-0.07 (0.18)\\
$\sigma_{23} =0$&  0.04 (0.22)& 0.03 (0.22)& 0.03 (0.25)\\
$\sigma_{33} =1$&  0.76 (0.72)& 0.88 (0.77)& 0.75 (0.76)\\
\hline
\end{tabular}
\end{center}
\end{table}

\begin{table}[p]
\caption{Bias (std) under 3-component trivariate
normal mixture models.}
\label{tab4.4b}
\begin{center}
\vspace{1ex}
\renewcommand{\arraystretch}{1.0}
\begin{tabular}{|l|c|c|c|}
  \hline
& {MLE} & {PMLE1} & {PMLE2} \\ \hline%
&\multicolumn{3}{c|}{Model IV.2.4, component 1} \\ \hline%
$\pi_1 =0.15   $&  0.00 (0.05)& 0.00 (0.04)& 0.01 (0.04)\\
$\mu_1    =0 $  &  0.04 (0.43)& 0.04 (0.37)& 0.02 (0.29)\\
$\mu_2    =0 $  &  0.20 (0.96)& 0.20 (0.90)& 0.24 (0.88)\\
$\mu_3    =-2$  &  0.19 (0.86)& 0.17 (0.80)& 0.20 (0.80)\\
$\sigma_{11} =1$&  0.05 (0.63)& 0.02 (0.52)& 0.01 (0.38)\\
$\sigma_{12} =0$& -0.03 (0.54)&-0.01 (0.41)&-0.01 (0.34)\\
$\sigma_{13} =0$&  0.04 (0.79)& 0.01 (0.58)& 0.01 (0.35)\\
$\sigma_{22} =1$&  0.18 (1.06)& 0.13 (0.81)& 0.18 (0.73)\\
$\sigma_{23} =0$& -0.15 (1.09)&-0.10 (0.65)&-0.09 (0.62)\\
$\sigma_{33} =1$&  0.65 (2.52)& 0.53 (2.17)& 0.68 (2.31)\\
\hline              %
&\multicolumn{3}{c|}{Model IV.2.4, component 2} \\ \hline%
$\pi_2 =0.35   $& -0.01 (0.06)&-0.01 (0.06)&-0.02 (0.06)\\
$\mu_1    =0$   &  0.01 (0.19)& 0.01 (0.19)& 0.01 (0.18)\\
$\mu_2    =3$   & -0.51 (1.25)&-0.46 (1.21)&-0.34 (1.13)\\
$\mu_3    =0$   &  0.24 (0.94)& 0.21 (0.91)& 0.13 (0.86)\\
$\sigma_{11} =1$&  0.56 (1.54)& 0.50 (1.47)& 0.35 (1.27)\\
$\sigma_{12} =0$& -0.49 (1.32)&-0.44 (1.26)&-0.32 (1.10)\\
$\sigma_{13} =0$&  0.09 (0.42)& 0.08 (0.42)& 0.05 (0.41)\\
$\sigma_{22} =3$&  0.48 (1.78)& 0.41 (1.71)& 0.20 (1.53)\\
$\sigma_{23} =0$& -0.33 (0.98)&-0.30 (0.96)&-0.25 (0.88)\\
$\sigma_{33} =10$&-1.40 (3.55)&-1.26 (3.45)&-1.03 (3.31)\\
\hline
&\multicolumn{3}{c|}{Model IV.2.4, component 3} \\ \hline%
$\pi_3 =0.5    $    &  0.01 (0.05)& 0.01 (0.05)& 0.00 (0.05)\\
$\mu_1    =0$       & -0.02 (0.18)&-0.02 (0.18)&-0.01 (0.19)\\
$\mu_2    =0$       &  0.37 (0.87)& 0.34 (0.86)& 0.27 (0.79)\\
$\mu_3    =2$       & -0.28 (0.72)&-0.25 (0.68)&-0.17 (0.58)\\
$\sigma_{11} =4.87$ & -0.57 (1.42)&-0.51 (1.36)&-0.39 (1.22)\\
$\sigma_{12} =-3.23$&  0.45 (1.24)& 0.41 (1.20)& 0.30 (1.07)\\
$\sigma_{13} =0.5$  & -0.07 (0.33)&-0.06 (0.33)&-0.04 (0.32)\\
$\sigma_{22} =7.2$  & -0.46 (1.48)&-0.42 (1.46)&-0.33 (1.38)\\
$\sigma_{23} =-2.16$&  0.31 (0.95)& 0.27 (0.89)& 0.18 (0.77)\\
$\sigma_{33}=1.94$  &  0.88 (2.23)& 0.79 (2.16)& 0.58 (1.90)\\
\hline
\end{tabular}
\end{center}
\end{table}

\begin{table}[p]
\caption{Bias (std) under 3-component trivariate
normal mixture models.}
\label{tab4.4c}
\begin{center}
\vspace{1ex}
\renewcommand{\arraystretch}{1.0}
\begin{tabular}{|l|c|c|c|}
  \hline
& {MLE} & {PMLE1} & {PMLE2} \\ \hline%
&\multicolumn{3}{c|}{Model IV.3.6, component 1} \\ \hline%
$\pi_1 =0.15    $    & 0.00  (0.05) & 0.00  (0.05) & 0.00  (0.05)\\
$\mu_1    =0 $       & 0.05  (0.41) & 0.05  (0.41) & 0.05  (0.40)\\
$\mu_2    =0 $       & -0.01 (0.64) & -0.01 (0.64) & -0.01 (0.61)\\
$\mu_3    =-2$       & -0.21 (1.23) & -0.21 (1.23) & -0.23 (1.20)\\
$\sigma_{11} =1 $    & 0.28  (1.24) & 0.28  (1.24) & 0.24  (1.12)\\
$\sigma_{12} =0 $    & -0.19 (1.16) & -0.19 (1.16) & -0.15 (1.05)\\
$\sigma_{13} =0 $    & 0.14  (1.04) & 0.14  (1.03) & 0.13  (0.99)\\
$\sigma_{22} =3 $    & 0.21  (1.48) & 0.21  (1.48) & 0.18  (1.40)\\
$\sigma_{23} =0 $    & -0.42 (1.54) & -0.42 (1.54) & -0.39 (1.50)\\
$\sigma_{33} =10$    & -1.37 (3.73) & -1.37 (3.73) & -1.34 (3.64)\\
\hline
&\multicolumn{3}{c|}{Model IV.3.6, component 2} \\ \hline%
$\pi_2 =0.35      $  & -0.01 (0.06) & -0.01 (0.06) & -0.01 (0.06)\\
$\mu_1    =0   $     & -0.01 (0.33) & -0.01 (0.33) & 0.00  (0.32)\\
$\mu_2    =3   $     & -0.20 (0.61) & -0.2  (0.61) & -0.19 (0.60)\\
$\mu_3    =0   $     & 0.25  (0.96) & 0.25  (0.96) & 0.26  (0.94)\\
$\sigma_{11}=4.87 $  & -0.15 (1.18) & -0.15 (1.18) & -0.13 (1.14)\\
$\sigma_{12}=-3.2 $  & 1.23  (2.89) & 1.23  (2.89) & 1.2   (2.87)\\
$\sigma_{13}=0.5  $  & -0.16 (0.62) & -0.16 (0.62) & -0.15 (0.62)\\
$\sigma_{22}=7.2  $  & -0.24 (1.56) & -0.24 (1.56) & -0.21 (1.52)\\
$\sigma_{23}=-2.16$  & 0.21  (0.77) & 0.21  (0.77) & 0.19  (0.73)\\
$\sigma_{33}=1.94 $  & 0.21  (1.61) & 0.21  (1.61) & 0.18  (1.52)\\
\hline              %
&\multicolumn{3}{c|}{Model IV.3.6, component 3} \\ \hline%
$\pi_3 =0.5        $ & 0.02  (0.07) & 0.02  (0.07) & 0.02  (0.07)\\
$\mu_1    =0    $    & -0.02 (0.22) & -0.02 (0.22) & -0.02 (0.22)\\
$\mu_2    =0    $    & 0.16  (0.43) & 0.17  (0.43) & 0.16  (0.43)\\
$\mu_3    =2    $    & -0.33 (0.68) & -0.33 (0.68) & -0.32 (0.68)\\
$\sigma_{11} = 4.87$ & -0.18 (0.66) & -0.18 (0.66) & -0.17 (0.65)\\
$\sigma_{12} = 3.23$ & -1.06 (2.14) & -1.06 (2.15) & -1.04 (2.15)\\
$\sigma_{13} =-0.5 $ & 0.17  (0.47) & 0.17  (0.47) & 0.16  (0.47)\\
$\sigma_{22} = 7.2 $ & -0.21 (0.97) & -0.21 (0.98) & -0.20 (0.98)\\
$\sigma_{23} =-2.16$ & 0.03  (0.45) & 0.03  (0.45) & 0.03  (0.46)\\
$\sigma_{33} =1.94 $ & 0.03  (0.39) & 0.03  (0.38) & 0.03  (0.38)\\
\hline
\end{tabular}
\end{center}
\end{table}

\end{document}